\documentclass[twoside]{article}

\usepackage{amsmath,amsthm,amssymb}

\usepackage{mathptmx}

\usepackage[text={12.5cm,19cm},centering,paperwidth=17cm,paperheight=24cm]{geometry}

\usepackage[fontsize=10.4pt]{scrextend}

\def\titlerunning#1{\gdef\titrun{#1}}

\makeatletter
\def\author#1{\gdef\autrun{\def\and{\unskip, }#1}\gdef\@author{#1}}
\def\address#1{{\def\and{\\\hspace*{15.6pt}}\renewcommand{\thefootnote}{}\footnote{#1}}\markboth{\autrun}{\titrun}}
\makeatother

\def\email#1{email: \href{mailto:#1}{#1} }

\newenvironment{dedication}{\itshape\center}{\par\medskip}
\newenvironment{acknowledgments}{\bigskip\small\noindent\textit{Acknowledgments.}}{\par}

\newtheorem{thm}{Theorem}[section]

\theoremstyle{definition}

\numberwithin{equation}{section}

\usepackage[hyperfootnotes=false,colorlinks=true,allcolors=blue]{hyperref}

\begin{document}

\titlerunning{Large $|k|$ behavior of CGO solutions}

\title{\textbf{Large $|k|$ behavior of  d-bar problems for domains with a
smooth boundary}}

\author{Christian Klein \and Johannes Sj\"ostrand \and Nikola Stoilov}

\date{}

\maketitle

\address{C. Klein: Institut de Math\'ematiques de Bourgogne, UMR 5584\\
                Universit\'e de Bourgogne-Franche-Comt\'e, 9 avenue Alain Savary, 21078 Dijon
                Cedex, France; \email{Christian.Klein@u-bourgogne.fr} 
				\and J. Sj\"ostrand: Institut de Math\'ematiques de Bourgogne, UMR 5584\\
                Universit\'e de Bourgogne-Franche-Comt\'e, 9 avenue Alain Savary, 21078 Dijon
                Cedex, France; \email{Johannes.Sjostrand@u-bourgogne.fr}
				\and N. Stoilov: Institut de Math\'ematiques de Bourgogne, UMR 5584\\
                Universit\'e de Bourgogne-Franche-Comt\'e, 9 avenue Alain Savary, 21078 Dijon
                Cedex, France; \email{Nikola.Stoilov@u-bourgogne.fr}}

\begin{dedication}
To Ari Laptev on the occasion of his 70th birthday
\end{dedication}

\begin{abstract}
	In a previous work on the large $|k|$ behavior of 
	complex geometric optics solutions to a system of d-bar equations, 
	we treated in detail the situation when a certain
        potential is the characteristic function of a strictly convex set with
      real-analytic boundary. We here extend the results to the case
      of sets with smooth boundary, by using almost holomorphic functions.  
\end{abstract}

%
%

\section{Introduction}

\par
This note is concerned with solutions to the  Dirac system 
\begin{equation}\label{dbarphi}
  \begin{cases}
    \bar{\partial}\phi_{1}=\frac{1}{2}q\mathrm{e}^{\bar{k}\bar{z}-kz}\phi_{2},\\
    \partial\phi_{2}=\sigma\frac{1}{2}\bar{q}\mathrm{e}^{kz-\bar{k}\bar{z}}\phi_{1},\quad \sigma=\pm1,
\end{cases}   
\end{equation}
subject to the asymptotic conditions 
\begin{equation}
    \lim_{|z|\to\infty}\phi_{1}=1,\quad \lim_{|z|\to\infty}\phi_{2}=0;
    \label{Phisasym}
\end{equation}
where  $q=q(x,y)$ is a complex-valued field, the \emph{spectral 
parameter} $k\in\mathbb{C}$ is independent of $z=x+\mathrm{i} y$, and
\begin{equation*}
\partial:=\frac{1}{2}\left(\frac{\partial}{\partial x}-\mathrm{i}\frac{\partial}{\partial y}\right)\quad\text{and}\quad
\bar{\partial}:=\frac{1}{2}\left(\frac{\partial}{\partial 
x}+\mathrm{i}\frac{\partial}{\partial y}\right).
\end{equation*}
The functions $\phi_{i}(z;k)$, $i=1,2$  depend on $z$ and $k$ where it 
is understood that they need not be holomorphic in either variable. 
This system appears in the scattering theory of the integrable 
Davey-Stewartson II equation \cite{Fok} and in the context of 
Calder\'on's problem \cite{calderon} (for more applications and 
references, see \cite{KlSjSt20}). 

In \cite{KlSjSt20}, we have addressed the large $|k|$ behavior of the 
solutions to (\ref{dbarphi}) and (\ref{Phisasym}). It was shown that 
the $k$ dependence for large $|k|$ strongly depends on the regularity 
of the potential $q$ in (\ref{dbarphi}). Of particular interest is 
the case of potentials being the characteristic function of a compact 
domain with strictly convex boundary. It was shown in \cite{KlSjSt20} 
that in this case $\phi_{2}=\mathcal{O}(|k|^{-1/2})$ in polynomially
weighted $L^2$-spaces. The leading 
order  contribution to $\phi_{2}$ is equal to
  $\overline{f(z,k)}/(2\pi) $, where $f$ is given by the integral 
\begin{equation}
f(z,k) = \int_\Omega \frac{1}{z - w}e^{\overline{kw} -kw}L(\mathrm{d}w) = \iint_\Omega \frac{e^{\overline{kw} -kw}}{z - w}\frac{\mathrm{d}\overline{w}\wedge\mathrm{d}w}{2i},
\end{equation}
for $z,k \in \mathbb{C}$, $|k| \gg 1$.  Here 
\begin{equation*}
L(\mathrm{d}w) = \mathrm{dRe}(w)\wedge \mathrm{dIm}(w) = \frac{\mathrm{d}\overline{w}\wedge\mathrm{d}w}{2i}
\end{equation*}
is the Lebesgue measure. The large $|k|$ asymptotics of $f$ was
computed in \cite{KlSjSt20} when $\partial \Omega $ is analytic.

In this note we extend the results of Section 5 in \cite{KlSjSt20} to
the case of strictly convex domains $\Omega \Subset \mathbb{C}$ with a
boundary that is smooth but not necessary analytic. In
  doing so, we follow the complex analysis arguments in
  \cite{KlSjSt20}, but make an essential use of almost holomorphic
  extensions, introduced by H\"ormander \cite{Ho69} and others. See
  \cite{Ho69} with the attached commentary and additional bibliography by one of
  us. Almost holomorphic extensions play an important role in the
  global theory of Fourier integral operators with complex phase
  \cite{MeSj75}. Complex phase arguments were also used by
  A.\ Laptev, Yu.\ Safarov, and D.\ Vassiliev \cite{LaSaVa94} in their
  approach to Fourier integral operators with real phase and their
  global representations as superpositions of Gaussian type operators
  and at one point the authors evoke almost analytic (i.e.\ almost
  holomorphic) extensions.
  
  \section{Main result}

As in Section 5.2. in  \cite{KlSjSt20} we use Stokes' formula to get 
\begin{equation}
f(z,k) =
\frac{1}{2i\overline{k}}\int_{\partial\Omega}\frac{1}{z-w}e^{\overline{kw}
  -kw}\mathrm{d}w + \frac{\pi}{\overline{k}}e^{\overline{kw}
  -kw}1_\Omega (z)
\label{ic.1}.
\end{equation}
We write $kw - \overline{kw} = (i/h)\Re(w\overline{\omega})$, with $h = 1/|k|$, $\omega =  2\overline{k}/|k| $.

Consider the function $iu_0(w) = kw - \overline{kw}$ on $\partial\Omega$ (suppressing sometimes the large parameter from the notation). In \cite{KlSjSt20} we use the analyticity of $\partial\Omega$ to introduce a holomorphic extension $iu(w)$ of $iu_0(w)$  from $\partial\Omega$ to $\mathrm{neigh}(\partial\Omega, \mathbb{C} )$.

When $\partial\Omega$ is only smooth we take instead an almost holomorphic extension $u\in C^\infty(\mathrm{neigh}(\partial\Omega, \mathbb{C} ))$ so that
\begin{equation}
u|_{\partial\Omega} = u_0, \quad \partial_{\overline{w}}u(w) = \mathcal{O}(\mathrm{dist}(w, \partial\Omega)^\infty)
\label{ic.101},
\end{equation} 
see \cite{MeSj75},\cite{SjZw91}. We recall that $u$ is unique up to infinite order: If $\tilde{u}$ is another analytic holomorphic extension of $u_0$ then
\begin{equation}
u - \tilde{u} = 
\mathcal{O}(\mathrm{dist}(\cdot,\partial\Omega)^\infty)~\text{in}~ \mathrm{neigh}(\partial\Omega, \mathbb{C}).
\end{equation}
 As in \cite{KlSjSt20} eq. (5.13) we parametrize points $w$ in a neighbourhood of $\partial\Omega$ by 
 \begin{equation}
 w = \gamma(t)+ \mathrm{i}s\dot{\gamma}(t), \quad t\in 
 \mathbb{R}/L\mathbb{Z}, \quad s \in\mathrm{neigh}(0,\mathbb{R}) ,
 \end{equation}
 where $L>0$ is the length of $\partial\Omega$. Here $t\to \gamma(t)$ 
 is a parametrization of $\partial\Omega$ with positive orientation 
 and with $|\dot{\gamma}(t)| = 1$, such that $\gamma$ can be viewed 
 as an $L$-periodic function. As in \cite{KlSjSt20} Section 5, we can 
 view $u$ as a function $u = u(t,s)$, $t\in \mathbb{R}/L\mathbb{Z}$, $s\in \mathrm{neigh}(0, \mathbb{R})$, satisfying 
 \begin{equation}
 \partial_t u + \mathrm{i}\left(1+ 
 \mathrm{i}s\frac{\ddot{\gamma}}{\gamma}\right)\partial_su = \mathcal{O}(s^\infty), 
\end{equation}  
and with $u(t,s) = u_0(t) + su_1(t)+s^2u_2(t)+ \mathcal{O}(s^3)$, we get

\begin{equation}\label{iche.15}
\partial _tu_0(t)=\frac{|k|}{2}\left(\dot{\gamma }(t)\overline{\omega }
    +\overline{\dot{\gamma }}(t) \omega  \right)=|k|\langle \dot{\gamma
  }(t),\omega \rangle_{\mathbb{R}^2},
\end{equation}
\begin{equation}\label{iche.16}
\partial^2 _tu_0(t)=\frac{|k|}{2}\left(\ddot{\gamma }(t)\overline{\omega }
    +\overline{\ddot{\gamma }}(t) \omega  \right)=|k|\langle \ddot{\gamma
  }(t),\omega \rangle_{\mathbb{R}^2}.
\end{equation}

\par From (\ref{iche.15}) we see that $\gamma (t)$ is a critical point
of $u_0(t)$ iff $\omega $ (which is non-vanishing) is normal to
$\partial \Omega $ at $\gamma (t)$. From $\langle \dot{\gamma }(t),
\dot{\gamma }(t)\rangle =1$ we know that
\begin{equation}\label{iche.17}
\langle \dot{\gamma }(t),\ddot{\gamma }(t)\rangle =0,
\end{equation}
and hence $\ddot{\gamma }(t)$ is normal to $\partial \Omega $
everywhere. Thus at a critical point of $u_0$ we have $\ddot{\gamma
}(t)\in \mathbb{R}\omega $. It follows from (\ref{iche.16}) that a critical
point is non-degenerate precisely when $\ddot{\gamma }(t)\ne 0$, i.e.\
when $\partial \Omega $ has non-vanishing curvature there. Such a
point is
\begin{itemize}
  \item a local maximum if $\ddot{\gamma }(t)=c\omega $, $c<0$, and
  \item a local minimum if $\ddot{\gamma }(t)=c\omega $, $c>0$.
\end{itemize}

\par Now recall the assumption that
\begin{equation}\label{iche.18}
\Omega \hbox{ is strictly convex}.
\end{equation}
Then at every point in $\partial \Omega $, $\ddot{\gamma }(t)$ is
non-vanishing and of the form $c(t)\nu (t)$, where $c(t)>0$ and
$\nu (t)=\mathrm{i}\dot{\gamma }(t)$ is the interior unit normal (recall that
$\gamma $ is positively oriented).

\par For a fixed $k\ne 0$, we can
decompose
\begin{equation}\label{iche.19}
\partial \Omega =\{ w_-(k )\}\cup \Gamma _+\cup \{ w_+(k )\}\cup
\Gamma _-,
\end{equation}
ordered in the positive direction when starting and ending at
$w_-(k )$. Here
\begin{itemize}
\item $w_-(k )$ is the south pole, where $\nu =c\omega $ for some
  $c<0$. Equivalently this is the global maximum point of $u_0$.
  \item $\Gamma _+$ is the open boundary segment connecting
    $w_-(k )$ to $w_+(k )$ in the positive direction.
\item $w_+(k )$ is the north pole, where $\nu =c\omega $ for some
  $c>0$. Equivalently this is the global minimum point of $u_0$.
  \item $\Gamma _-$ is the open boundary segment connecting
    $w_+(k )$ to $w_-(k )$ in the positive direction.
  \end{itemize}
  Notice that
  \begin{equation}\label{iche.20}
\Gamma _\pm =\{ \gamma (t)\in \partial \Omega ;\, \mp \partial
_tu_0(\gamma (t))>0 \} .
\end{equation}

\medskip
\par On $\partial \Omega $ we have :
$$
e^{\overline{kw}-kw}=e^{-iu_0(w)},\ u_0(w)=|k|\Re
(w\overline{\omega  }),\ w=\gamma (t).
$$
The formula (\ref{ic.1}) reads
\begin{equation}\label{icbtf.1}
f(z,k)=\frac{1}{2i\overline{k}}\int_{\partial \Omega
}\frac{1}{z-w}e^{-iu(w,k )}dw+
  (\pi/\overline{k})e^{-i|k|\Re (z\overline{\omega })}1_\Omega(z) ,
\end{equation}
where $u(\cdot ,k)$ is an almost holomorphic extension of
$u_0(w)=u_0(w,k)$ to a neighbourhood of $\partial \Omega $.

We introduce the same closed contour $\Gamma$ as in \cite{KlSjSt20}, 
Subsection 5.3. Recall that $\Gamma$ is a deformation of 
$\partial\Omega$ obtained by pushing $\Gamma_+$ inward and $\Gamma_-$ 
outward. This is obtained by means of a smooth vector
  field $\upsilon $, defined near $\partial \Omega $, transversal to
  $\partial \Omega $ and pointing outward, so that 
\begin{equation*}
\Gamma = \{ \exp (t\upsilon (w))|_{t = \tau(w)}; w\in\partial\Omega\},
\end{equation*} 
when $\tau(w)$ is a suitable function such that $\sup|\tau|\ll 1$, $\tau<0$ on $\Gamma_+$, $\tau>0$ on $\Gamma_-$. 
Near $w_+$ we can choose a real Morse coordinate $\mu$ on $\Gamma$, centred at $w_+$ such that as a function of $\mu$,
\begin{equation}
u = u(w_+(k), k) + |k|\mu^2/2.
\label{ichtf.7}
\end{equation}
Let $\mu$ also denote an almost holomorphic extension. Then in a complex neighbourhood of $w_+$, we have $u$ as a function of $\mu$
\begin{equation}
u = u(w_+ (k),k) + |k|\mu^2/2 +|k|\mathcal{O}((\Im\mu)^\infty).
\label{ichtf.71}
\end{equation}
Near $w_+$ we choose $\Gamma$ to coincide with
$\exp(-\mathrm{i}\pi/4)\mathbb{R}$  in the complex $\mu $-coordinate. Thus along this part of $\Gamma$, we have 
\begin{equation*}
-\mathrm{i}\left(u(\cdot ,k) - u(w_+(k),k)\right) = - |k||\mu|^2/2 +
|k|\mathcal{O}(|\mu |^\infty).
\end{equation*}
Near $w_-$ we have a similar construction. Along $\Gamma$ we have 
\begin{equation}
\label{icbtf.8}
|e^{-iu(w,k )}|\le e^{-|k |\mathrm{dist\,}(w,\{w_+(k ),w_-(k)\} )^2/C}.
\end{equation} 

Let
$\Omega _+\subset \overline{\Omega } $ and $\Omega _-\subset \mathbb{C}\setminus \Omega $, be the points swept over by the deformation of
$\Gamma _+$ and $\Gamma _-$ respectively,
$$\Omega _+=\{ \exp t\upsilon (w);\, w\in \Gamma _+,\ \tau (w)\le t
\le 0  \},$$
$$\Omega _-=\{ \exp t\upsilon (w);\, w\in \Gamma _-,\ 0\le t
\le \tau (w)  \},$$
and notice that
$$
\Gamma =\partial ((\Omega \setminus \Omega _+)\cup \Omega _-).
$$

Assume for simplicity that $z\not\in \partial \Omega_{+} \cup \partial
\Omega_{-} $ and put
\begin{equation}\label{FGam}
F(z)=F_\Gamma (z)=\int_\Gamma \frac{1}{z-w}e^{-iu(w,k )}dw.
\end{equation}
When $\Gamma$ is real-analytic, we can choose $u$ holomorphic and using the residue theorem we got in \cite{KlSjSt20}
\begin{equation}\label{wn.3}
  \begin{split}
f(z, & k)=\frac{1}{2 i \overline{k}}\int_{\partial \Omega
}\frac{1}{z-w}e^{-iu(w,k )}dw+(\pi /\overline{k})e^{-i|k|\Re
  (z\overline{\omega })}1_\Omega (z)\\
&=\frac{1}{2i\overline{k}}F(z)+(\pi
/\overline{k})\left(e^{-iu(z,k )}(1_{\Omega _-}(z)-1_{\Omega
    _+}(z))+e^{-i|k|\Re (z\overline{\omega })}1_\Omega (z) \right).
  \end{split}
\end{equation}
In the smooth case we get a correction term that we compute with the help of Stokes' theorem. Let $V\Subset \mathbb{C}$ have a smooth positively oriented boundary. If $g$ is a distribution defined near $\overline{V}$ sufficiently smooth near $\partial V$, Stokes' formula gives 
\begin{equation*}
\frac{1}{2\mathrm{i}}\int_{\partial V} g(w)\mathrm{d}w = \int_V 
\partial_{\overline{w}}g \frac{\mathrm{d}\overline{w}\wedge\mathrm{d}w}{2\mathrm{i}} = \int_{V}\partial_{\overline{w}}gL(\mathrm{d}w).
\end{equation*}
With $g(w) = \frac{1}{z-w}e^{-\mathrm{i}u(w,k)}$  and
  $\delta _z$ denoting the delta mass at the point $z$, we get
\begin{equation*}
\partial_{\overline{w}}g = -\pi\delta_z(w)e^{-\mathrm{i}u(z,k)} - \frac{\mathrm{i}\partial_{\overline{w}}u}{z-w}e^{-iu(w,k)}, 
\end{equation*}
hence
\begin{equation*}
\frac{1}{2\mathrm{i}}\int_{\partial V}\frac{1}{z-w}e^{-\mathrm{i}u(w,k)\mathrm{d}w} = -\pi e^{-iu(z,k)} 1_V(z) - \mathrm{i}\int_V\frac{\partial_{\overline{w}}u(w,k)e^{-\mathrm{i}u(w,k)}}{z-w}L(\mathrm{d}w)
\end{equation*}
We apply this with $V$ equal to $\Omega_+$ and $\Omega_-$ and get 
\begin{multline*}
\frac{1}{2\mathrm{i}}\int_{\partial\Omega}\frac{1}{z-w}e^{-\mathrm{i}u(w,k)}\mathrm{d}w - \frac{1}{2\mathrm{i}}\int_\Gamma\frac{1}{z-w}e^{-\mathrm{i}u(w,k)}\mathrm{d}w = \\
 -\pi e^{-iu(z,k)}1_{\Omega_+}(z) + \pi e^{-\mathrm{i}u(z,k)}1_{\Omega_-}(z) 
 -\mathrm{i}\int\frac{\partial_{\overline{w}}u(w,k)}{z-w}e^{-\mathrm{i}u(w,k)}
(1_{\Omega _+}(w)-1_{\Omega _-}(w))
 L(\mathrm{d}w) 
\end{multline*}
The last term gives rise to a new term in (\ref{wn.3}) and we get
\begin{multline}\label{wn.4}
  f(z,k) = \frac{1}{2i\overline{k}}F(z)+(\pi
  /\overline{k})\left(e^{-iu(z,k )}(1_{\Omega _-}(z)-1_{\Omega
      _+}(z))+e^{-i|k|\Re (z\overline{\omega })}1_\Omega (z)
  \right)\\
  +
  \frac{i}{\overline{k}}\int\frac{\partial_{\overline{w}}u(w,k)}{z-w}e^{-\mathrm{i}u(w,k)}\left(1_{\Omega_-}(w)
    - 1_{\Omega_+}(w)\right)L(\mathrm{d}w). \end{multline} In order to
estimate the last term we notice that
\begin{equation}
\left|e^{-\mathrm{i}u(w,k)}\right|\leq e^{-|k|\mathrm{dist}(w, \partial\Omega)/C}, \quad w\in(\Omega_+\cup\Omega_-)\setminus\mathrm{neigh}(\{ w_+, w_-\}).
\label{wn.5}
\end{equation}
In a neighbourhood of $w_+$, we use the Morse coordinates in (\ref{ichtf.71}) and notice that with $t = \Re(\mu)$, $s = \Im(\mu)$ we have in 
\begin{equation*}
\mathrm{neigh}(w_+)\cap\Omega_-:0<t<1/\mathcal{O}(1),\quad -t<s<0,
\end{equation*}
and in 
\begin{equation*}
\mathrm{neigh}(w_+)\cap\Omega_+:-1/\mathcal{O}(1)<t<0,\quad 0<s<-t.
\end{equation*}

In both sets, $-\Im(u)  = - ts +\mathcal{O}(s^\infty)\geq s^2-\mathcal{O}(s^\infty)$, so
\begin{equation}
|e^{-\mathrm{i}u(w,k)}|\leq e^{-|k|(s^2+\mathcal{O}(s^\infty))} \text{ in } \mathrm{neigh}(w_+)\cap(\Omega_-\cup\Omega_+) 
\label{wn.6}.
 \end{equation}
The same estimate holds  near $w_-$ in  the corresponding Morse 
coordinate. Combining (\ref{wn.5}), (\ref{wn.6}), (\ref{ic.101}) we get 
\begin{equation*}
|(\overline{\partial}_w u ) e^{\mathrm{i}u(w,k)}| =
\mathcal{O}(|k|^{-\infty})~ \text{ uniformly in }~
\Omega_+\cup\Omega_- .
\end{equation*} 
We conclude that the last term in (\ref{wn.4}) is
$\mathcal{O}(<z>^{-1}|k|^{-\infty})$.

\par The analysis of $F_\Gamma $ in \cite{KlSjSt20} goes through
without any changes and we obtain:
\begin{thm}\label{Th_smooth}
The detailed description of the asymptotics of $f$ when $k\to \infty $
in \cite[Theorem 5.2]{KlSjSt20} remains valid if we replace the
assumption that $\partial \Omega $ is real-analytic, by the weaker
assumption that $\partial \Omega $ is smooth.
\end{thm}
The results in \cite[Subsection 5.5]{KlSjSt20} about the
leading approximations of $\phi _1$, $\phi _2$ in certain polynomially
weighted $L^2$-norms remain valid after the same weakening of the
boundary regularity assumption for $\Omega $.

\begin{acknowledgments}
	This work is partially supported by 
the ANR-FWF project ANuI - ANR-17-CE40-0035, the isite BFC project 
NAANoD, the EIPHI Graduate School (contract ANR-17-EURE-0002) and by the 
European Union Horizon 2020 research and innovation program under the 
Marie Sklodowska-Curie RISE 2017 grant agreement no. 778010 IPaDEGAN.
\end{acknowledgments}

\end{document}